\documentclass[11pt]{article}
\usepackage{amssymb}
\newtheorem{conj}{Conjecture}[section]
\newtheorem{thm}[conj]{Theorem}
\newtheorem{rem}[conj]{Remark}

\begin{document}
\title{\Large On Tverberg's conjecture}

\author{{\Large Sini\v sa T. Vre\' cica}
\footnote{Supported by the Ministry for Science,
Technology and Development of Serbia, Grant 1854}}
\maketitle

\begin{abstract}
In 1989 Helge Tverberg proposed a quite general conjecture in
Discrete geometry, which could be considered as the common basis
for many results in Combinatorial geometry and at the same time as
a discrete analogue of the common transversal theorems. It implies
or contains as the special cases many classical "coincidence"
results such as Radon's theorem, Rado's theorem, the Ham-sandwich
theorem, the "nonembeddability" results (e.g. nonembeddability of
graphs $K_5$ and $K_{3,3}$ in $\mathbb R^2$) etc.

The main goal of this short note is to verify this conjecture in
one new, non-trivial case. We obtain the continuous version of the
conjecture. So, it is not surprising that we use the topological
methods, or more precisely the methods of equivariant topology and
the theory of characteristic classes.
\end{abstract}

\section{Introduction}

Establishing the relation between Rado's theorem on general measure (see \cite{r})
and the Ham sandwich theorem, the following result is proved in \cite{zv} stating
that these two results belong to the same family.

\begin{thm}
\label{ZV}
Let $0\leq k\leq d-1$, and let $\mu_0,\mu_1,...,\mu_k$ be
$\sigma$-additive probability measures on $\mathbb R^d.$ Then there is a
$k$-flat $F$ with the property that every
closed halfspace containing $F,$ has $\mu_i$-measure at least \, $\frac
1{d-k+1}$ \, for all \, $i, \, 0\leq i\leq k.$
\end{thm}

Namely, this theorem reduces to Rado's theorem in the case $k=0$ and to the Ham
sandwich theorem in the case $k=d-1$.

The proof of the above theorem uses the topological result claiming the
non-existence of a non-zero section of a certain vector bundle over the Grassmann
manifold. Helge Tverberg observed that the special case $k=0$ (Rado's theorem)
follows easily from his result in combinatorial geometry from \cite{t}.

\begin{thm}
\label{T}
Let $S$ be a set of $(r-1)(d+1)+1$ points in $\mathbb R^d.$ Then one can split
it into subsets \, $S_1,S_2,...,S_r$ \, so that
$$\cap_{i=1}^r \mbox { conv } S_i\neq \emptyset.$$
\end{thm}

This observation motivated him to suppose that a general result should exist
which would generalize his result \ref{T} and at the same time  it would imply
theorem \ref{ZV} in the same way as his result implied Rado's theorem. So, he
formulated the following {\bf Tverberg's conjecture}:

\begin{conj}
\label{Tc}
Let $0\leq k\leq d-1$ and let $S_0,S_1,...,S_k$
be finite sets of points in $\mathbb R^d,$ with \, $\mid
S_i\mid=(r_i-1)(d-k+1)+1$ \, for \, $i=0,1,...,k.$ Then $S_i$ can be
split into $r_i$ sets, \, $S_{i1},S_{i2},...,S_{ir_i},$ so that there
is a $k$-flat $F$ meeting all the sets $\mbox { conv } S_{ij},
\, 0\leq i\leq k, \, 1\leq j\leq r_i.$
\end{conj}

It is easy to see that this conjecture implies the theorem
\ref{ZV} (\cite{t2}, see also \cite{z}) and its special case $k=0$
is the theorem \ref{T}. This conjecture unifies two important
themes of Combinatorial geometry: Helly-type theorems (the special
case $r=2$ of the theorem \ref{T} is the well known Radon's
theorem), and the common transversal theorems. Moreover, it was
shown in \cite{z} that this result could be considered as an
example of the whole family of results of the "combinatorial
geometry on vector bundles", and that these results would also
generalize many coincidence results such as the nonembeddability
of the graphs having minor $K_5$ or $K_{3,3}$ in the plane etc.

In paper \cite{tv} this conjecture was verified in some special
cases. It was proved that, besides the already known case $k=0$
the conjecture is true in the case $k=d-1$, and in the case $k=1$
when $r_0=1$ or $r_1=1$ or $r_0=r_1=2$. Also, slightly weakened
version of the conjecture is proved in the case $k=d-2$, obtained
when $3r_i$ points in $S_i$ were considered instead of $3r_i-2$ of
them. This version still suffices to imply theorem \ref{ZV} in the
case $k=d-2$.

Much more general result verifying some cases of this conjecture was achieved by
R. \v Zivaljevi\' c in \cite{z}, where he established the conjecture \ref{Tc} when
$r_0=r_1=\cdots =r_k$ is an odd prime number, and both $d$ and $k$ are odd integers.

The main result of this paper is the establishing of the conjecture \ref{Tc}
in the case $r_0=r_1=\cdots =r_k=2$ without any additional restriction.

The result in \cite{z} is obtained by using parametrized, ideal-valued,
cohomological index theory. The method of this paper will again be the reduction
of the result to the non-existence of some equivariant mapping, or equivalently
to the non-existence of the section of certain sphere bundle. This
will be proved by showing that the corresponding Stiefel-Whitney class is
non-trivial.

\section{The result}

As we already mentioned, the main goal of this paper is to
establish the conjecture \ref{Tc} in the case $r_0=r_1=\cdots
=r_k=2$, i.e. to prove the following theorem.

\begin{thm}
\label{linear}
Let $0\leq k\leq d-1$ and let $S_0,S_1,...,S_k$
be finite sets each of $d-k+2$ points in $\mathbb R^d$.
Then every set $S_i$ can be
split into $2$ sets, \, $S_{i1}$ and $S_{i2},$ so that there
is a $k$-flat $F$ meeting all the sets $\mbox { conv } S_{ij},
\, 0\leq i\leq k, \, j\in \{1,2\}.$
\end{thm}

Moreover, we will prove a slight generalization of this result. Namely, we could
consider the family of \, $(d-k+1)$-dimensional simplicies $\Delta_0,
\Delta_1,...,\Delta_k$, where $\Delta_i$ is the simplex spanned by the vertices
$e^i_0,e^i_1,...,e^i_{d-k+1}$ for $i=0,1,...,k$. We consider then the linear
mapping sending the vertices of $\Delta_i$ to the points in $S_i$ and note
that the sets $\mbox { conv } S_{i1}$ and $\mbox { conv } S_{i2}$ are actually
the images of the pair of disjoint faces of the simplex $\Delta_i$ under the
considered linear mapping. We observe that our argument also works for any
continuous mapping and not only for linear ones. So, we prove the following theorem.

\begin{thm}
\label{main}
Let $0\leq k\leq d-1$ and let $\varphi_i : \Delta_i \rightarrow \mathbb R^d$,
for $i=0,1,...,k$ be continuous mappings. Then there is an affine $k$-flat $F$
which intersects the images of two disjoint faces of every simplex $\Delta_i$.
\end{thm}

\medskip\noindent
{\bf Proof:}
Let us denote with $\mbox{Gr}_{d,d-k}$ the Grassmann manifold of
$(d-k)$-dimensional linear subspaces of $\mathbb R^d$, and for any
$L\in \mbox{Gr}_{d,d-k}$, with ${\pi}_L : \mathbb R^d \rightarrow L$ the
orthogonal projection.

Let us consider $(d-k+1)$-dimensional simplicies $\Delta_i=\mbox{
conv } \{e^i_0,e^i_1,...,e^i_{d-k+1}\}, \; i=0,1,...,k$. With
$(\Delta_i)^2_{\ast}$ we denote the deleted square of the simplex
$\Delta_i$, i.e. the set of ordered pairs of points in $\Delta_i$
having disjoint supports. It is easy to verify that the deleted
square $(\Delta_i)^2_{\ast}$ is a $(d-k)$-dimensional manifold
which we denote by $M$. It could be proved that $M$ is actually
homeomorphic to the sphere $S^{d-k}$, but we don't need this fact
here. (It is shown in \cite{bss} that $M$ is $(d-k-1)$-connected.)

The mapping $\varphi_i : \Delta_i \rightarrow \mathbb R^d$ induces the mapping
$$\tilde{\varphi_i} : (\Delta_i)^2_{\ast} \rightarrow {\mathbb R^d}\times
{\mathbb R^d} \; , \; \tilde{\varphi_i}(x,y)=(\varphi_i(x),\varphi_i(y)).$$

Now we consider the product of these induced mappings
$\tilde{\varphi_i}$ and get the mapping
\smallskip

\begin{center}
$ \tilde{\varphi} :
(\Delta_0)^2_{\ast} \times \cdots \times (\Delta_k)^2_{\ast}
\rightarrow ({\mathbb R^d})^2 \times \cdots \times ({\mathbb
R^d})^2,$\\
\smallskip

$\tilde{\varphi}((x_0,y_0),...,(x_k,y_k))=(\tilde{\varphi_0}(x_0,y_0),
...,\tilde{\varphi_k}(x_k,y_k)).$
\end{center}
\smallskip

The statement of the theorem reduces to the claim that there
exists $L\in \mbox{Gr}_{d,d-k}$ such that $$
{\pi}_L(\varphi_0(x_0))={\pi}_L(\varphi_0(y_0))=\cdots
={\pi}_L(\varphi_k(x_k))={\pi}_L(\varphi_k(y_k)),$$

\noindent for some \, $((x_0,y_0),...,(x_k,y_k))\in
(\Delta_0)^2_{\ast}\times \cdots \times (\Delta_k)^2_{\ast}$. Here
$L=F^{\perp}$, i.e. $L$ is the orthogonal complement to the affine
$k$-flat $F$ claimed to exist in the statement of the theorem.

Let us denote with $\xi$ the canonical vector bundle over
$\mbox{Gr}_{d,d-k}$. For every $L\in \mbox{Gr}_{d,d-k}$, we have
the mapping
\smallskip

\begin{center}
$\psi_L=\pi_L\circ \tilde{\varphi} :
(\Delta_0)^2_{\ast} \times \cdots \times (\Delta_k)^2_{\ast}
\rightarrow L^{2k+2},$\\
\smallskip

$\psi_L((x_0,y_0),...,(x_k,y_k))=(\pi_L(\varphi_0(x_0)),\pi_L(\varphi_0(y_0))
,...,\pi_L(\varphi_k(x_k)),\pi_L(\varphi_k(y_k))).$
\end{center}
\smallskip

The group $G=\underbrace{\mathbb Z/2\oplus \cdots \mathbb
Z/2}_{k+1}$ acts on these spaces, freely on $(\Delta_0)^2_{\ast}
\times \cdots \times (\Delta_k)^2_{\ast}\approx M^{k+1}$
and fiberwise on ${\xi}^{2k+2}$ (i.e. trivially on
$\mbox{Gr}_{d,d-k}$). The above mapping is equivariant and it
induces a section $s$ of the vector bundle $$ {\xi}^{2k+2}\times_G
M^{k+1} \rightarrow \left(\mbox{Gr}_{d,d-k}\times
M^{k+1}\right)/G= \mbox{Gr}_{d,d-k}\times M^{k+1}/G.$$

The fiber over $[L,(x_0,y_0,...,x_k,y_k)]$ could be identified
with $L^{2k+2}.$ The statement of the theorem reduces now to the
claim that the section $s$ intersects the diagonal $\Delta$ in
some fiber. Let us suppose, to the contrary, that $s$ does not
intersect the diagonal in any fiber. Projecting to the orthogonal
complement of the diagonal and then radially to its sphere (in
each fiber), we obtain the non-zero section of the vector bundle
with the fiber $L^{2k+1}$ and the section of the associated sphere
bundle whose fiber is homeomorphic with $S^{(d-k)(2k+1)-1}.$

We reach a contradiction (proving in this way the theorem) by
showing that the top-dimensional Stiefel-Whitney class of this
sphere bundle does not vanish. The Poincare dual of the
top-dimensional Stiefel-Whitney class of the sphere bundle
coincides with the homology class of the zero-set of a section of
the associated vector bundle (with the fiber $L^{2k+1}$), which is
transversal to the zero section. So, it suffices to find a section
of the corresponding vector bundle which intersects the zero
section transversally in an odd number of points.

In order to construct such a section we consider $k+1$ parallel
$(d-k)$-dimensional affine planes $A_0,A_1,...,A_k$ in $\mathbb R^d$,
and $d-k+2$ points in each of them being the vertices
$v_0^i,v_1^i,...,v_{d-k}^i$ and a barycenter $\hat{\sigma_i}$ of
$(d-k)$-dimensional simplex $\sigma_i,\,  0\leq i\leq k$. Let us
suppose that they are in a generic position meaning that their
$k+1$ barycenters span an affine $k$-dimensional flat $F$. We also
consider linear mappings $f_i$ which map vertices of $\Delta_i$ to
the vertices and the barycenter of $\sigma_i,$ namely
$f_i(e_j^i)=v_j^i, 0\leq j\leq d-k$ and
$f_i(e_{d-k+1}^i)=\hat{\sigma_i}$ for $i=0,1,...,k$.

The images under the linear mapping $f_i$ of two disjoint faces of
the simplex $\Delta_i$ are the convex hulls of the corresponding
vertices of the simplex $\sigma_i$ and its barycenter
$\hat{\sigma_i}$. The only nonempty intersection of the convex
hulls of two disjoint subsets of
$\{v_0^i,v_1^i,...,v_{d-k}^i,\hat{\sigma_i}\}$ is $$ \mbox{conv}
\{v_0^i,v_1^i,...,v_{d-k}^i\}\cap \mbox{conv}
\{\hat{\sigma_i}\}=\{\hat{\sigma_i}\}.$$

If some affine $k$-dimensional plane in $\mathbb R^d$ intersects
the images of two disjoint faces of some simplex $\Delta_i$, then
this plane contains the barycenter $\hat{\sigma_i}$ or it
intersects $A_i$ in at least $1$-dimensional affine plane. If this
plane should intersect the images of two disjoint faces of each
simplex $\Delta_0,\Delta_1,...,\Delta_k$, then (because of the
generic position of the planes $A_0,A_1,...,A_k$) it has to
contain the barycenters
$\hat{\sigma_0},\hat{\sigma_1},...,\hat{\sigma_k}$ and it is
uniquely determined by them. So, the mappings $f_0,f_1,...,f_k$
induce the section of the considered vector bundle which
intersects the diagonal at a single orbit of the action of the group
$G$, i.e. at the orbit $$ \left[L,\left(\frac {e_0^0+\cdots +e_{d-k}^0}
{d-k+1},e_{d-k+1}^0,...,\frac {e_0^k+\cdots +e_{d-k}^k}
{d-k+1},e_{d-k+1}^k\right)\right],$$

\noindent where $L=F^{\perp}$ is the orthogonal complement to the
$k$-dimensional flat $F$ spanned by the points
$\hat{\sigma_0},\hat{\sigma_1},...,\hat{\sigma_k}$.

So, we found a section which intersects the zero section of the
considered vector bundle in a single point. Obviously, for small
perturbation of the mappings $f_i$, the intersection with the zero
section will remain a single point, and so the obtained section
intersects the zero section transversally. This completes our
proof. \hfill $\square $

\begin{rem}
{\rm We note that we proved continuous version of the conjecture, and
so our theorem \ref{main} also generalizes the result of Bajm\' oczy
and B\'ar\' any from \cite{bb} in the same way as the result from
\cite{z} generalizes the result of B\' ar\' any, Shlosman and Sz\"
ucs from \cite{bss}.}
\end{rem}

%\nocite{*}
%\bibliographystyle{plain}
%\bibliography{t2}

\vskip 1cm

Sini\v sa T. Vre\' cica \par Faculty of Mathematics
\par University of Belgrade
\par Studentski trg 16, P.O.B. 550
\par 11000 Belgrade, Serbia
\par vrecica@matf.bg.ac.yu

\end{document}